\documentclass[a4paper,11pt]{article}
\usepackage{amsfonts,amsthm}
\usepackage[scale=0.7,centering]{geometry}

\newcommand{\E}{\mathbb{E}}
\newcommand{\erre}{\mathbb{R}}
\renewcommand{\H}{\mathcal{H}}
\newcommand{\K}{\mathcal{K}}

\newcommand{\ip}[2]{\left\langle#1,#2\right\rangle}
\newcommand{\ds}{\displaystyle}
\newcommand{\tr}{\mathop{\mathrm{Tr}}\nolimits}

\newtheorem{prop}{Proposition}[section]
\newtheorem{thm}[prop]{Theorem}
\newtheorem{lemma}[prop]{Lemma}
\newtheorem{coroll}[prop]{Corollary}
\newtheorem{defi}[prop]{Definition}
\theoremstyle{definition}
\newtheorem{rmk}[prop]{Remark}

\title{Variational inequalities in Hilbert spaces with measures and
  optimal stopping problems} 

\author{Viorel Barbu\footnote{University Al.~I.~Cuza, 8 Blvd. Carol
    I, Ia\c{s}i 700506, Romania. E-mail \texttt{vb41@uaic.ro}.}{\ }
  and Carlo Marinelli\footnote{Institut f\"ur Angewandte Mathematik,
    Universit\"at Bonn, Wegelerstr. 6, D-53115 Bonn, Germany. Tel.
    +49-178-918-3604, Fax +1-720-559-9266, e-mail
    \texttt{cm788@uni-bonn.de}. (Corresponding author)}}

\date{June 28, 2007}

\begin{document}

\maketitle

\begin{abstract}
  We study the existence theory for parabolic variational
  inequalities in weighted $L^2$ spaces with respect to excessive
  measures associated with a transition semigroup. We characterize the
  value function of optimal stopping problems for finite and infinite
  dimensional diffusions as a generalized solution of such a
  variational inequality. The weighted $L^2$ setting allows us to
  cover some singular cases, such as optimal stopping for stochastic
  equations with degenerate diffusion coefficient.  As an application
  of the theory, we consider the pricing of American-style contingent
  claims.  Among others, we treat the cases of assets with stochastic
  volatility and with path-dependent payoffs.

\medskip


\smallskip

\end{abstract}

\section{Introduction}
The aim of this work is to study a general class of parabolic
variational inequalities in Hilbert spaces with suitably chosen
reference measures. In particular, our motivation comes from the
connection between American option pricing in mathematical finance and
variational inequalities. It is well known by the classical works of
Bensoussan \cite{Ben} and Karatzas \cite{Kamerican} that the price of
an American contingent claim is the solution of an optimal stopping
problem, whose value function can be determined, in many cases,
solving an associated variational inequality (see e.g. \cite{GrSh} for
the classical theory and \cite{JLL} for connections with American
options).

In this paper we study variational inequalities associated to finite
and infinite dimensional diffusion processes in $L^2$ spaces with
respect to suitably chosen measures. In particular, denoting by $L$
the Kolmogorov operator associated to a diffusion $X$ on a Hilbert
space $H$, we shall choose a probability measure $\mu$ that is
(infinitesimally) excessive for $L$, i.e. that satisfies
$L^*\mu\leq\omega\mu$ for some $\omega\in\erre$ (see below for precise
statements). An appropriate choice of reference measure is essential
in the infinite dimensional case, as there is no analog of the
Lebesgue measure, and turns out to be useful also in the finite
dimensional case to overcome certain limitations of the classical
theory. In particular, we can relax the usual nondegeneracy
assumptions on the diffusion coefficient (or on the volatility, using
the language of mathematical finance), which is usually assumed in the
``traditional'' approach of studying variational inequalities in
Sobolev spaces w.r.t.  Lebesgue measure (see \cite{BL-ineq},
\cite{JLL}).  This allows us, for instance, to characterize the price
of American contingent claims on assets with degenerate or stochastic
volatility as the solution of a variational inequality. Similarly, we
can treat path-dependent derivatives, as well as claims on assets with
certain non-Markovian price evolutions, using the infinite dimensional
theory.  We would like to mention that Zabczyk \cite{Zinc} already
considered variational inequalities (called there Bellman inclusions)
in weighted spaces with respect to excessive measures, including
specific formulas for excessive measures and applications to American
option pricing. However, some of our results on existence of solutions
for the associated variatonal inequalities are more general (our
assumptions on the payoff function are weaker, we allow time-dependent
payoffs), and we explicitly construct a reference excessive measure in
many cases of interest. Let us also recall that a study of diffusion
operators in $L^p$ spaces with respect to invariant measures (i.e.
measures $\mu$ such that $L^*\mu=0$) has been initiated in
\cite{R-Lp}.

The main tool we rely on to study the above mentioned optimal stopping
problems is the general theory of maximal monotone operators in
Hilbert spaces.  However, we need some extensions of the classical
results, which are developed below and seem to be new. In particular,
we establish abstract existence results for variational inequalities
associated to the Kolmogorov operator of finite and infinite
dimensional diffusions (on these lines see also \cite{Barbu-OSNS} and
\cite{men-impulse}).

Variational inequalities connected to optimal stopping problems in
finance have also been studied in the framework of viscosity
solutions, see e.g. \cite{OkRe}, \cite{GaSw}. In particular in the
latter paper the authors consider the problem of optimal stopping in
Hilbert space and as an application they price American interest rate
contingent claims in the Go{\l}dys-Musiela-Sondermann model. Using the
approach of maximal monotone operators, at the expense of imposing
only very mild additional assumptions on the payoff functions, we are
able to obtain more regular solutions, which also have the attractive
feature of being the limit of iterative schemes that can be
implemented numerically. Moreover, the additional conditions on the
payoff function we need are satisfied in essentially all situations of
interest in option pricing.

The paper is organized as follows: in section \ref{sec:ex} we prove
two general existence results for the obstacle problem in Hilbert
spaces. In section \ref{sec:viosp} we relate these results with the
optimal stopping problem in Hilbert space. Applications to the pricing
of American contingent claims are given in section 4.

\section{Abstract existence results}
\label{sec:ex}
Let us first introduce some notation and definitions. Given any
Hilbert space $E$, we shall always denote by $|\cdot|_E$ its norm and
by $\ip{\cdot}{\cdot}_E$ its scalar product. Moreover, we define
$C([0,T],E)$ as the space of $E$-valued continuous functions on
$[0,T]$, and $W^{1,p}([0,T],E)$, $1\leq p \leq \infty$, as the space
of absolutely continuous functions $\varphi:[0,T]\to E$ with
$\frac{d\varphi}{dt}\in L^p([0,T],E)$.  The space of Schwarz'
distributions on a domain $\Xi \subset \erre^n$ will be denote by
$\mathcal{D}'(\Xi)$. Similarly, $W^{s,p}(\Xi)$ stands for the set of
functions $\phi:\Xi\to\erre$ that are in $L^p(\Xi)$ together with
their (distributional) derivatives of order up to $s$. Finally,
$\phi\in W^{s,p}_{loc}(\Xi)$ if $\phi\zeta\in W^{s,p}$ for all $\zeta
\in C^\infty_c(\Xi)$, the space of infinitely differentiable functions
on $\Xi$ with compact support.

Let $H$ be a Hilbert space and $\mu$ be a probability measure on $H$.
Denote by $\H$ the Hilbert space $L^2(H,\mu)$.  Let $(P_t)_{t\geq 0}$
be a strongly continuous semigroup on $\H$ with infinitesimal
generator $-N$. We shall assume that
$$
|P_t\phi|_\H \leq e^{\omega t} |\phi|_\H
\qquad \forall t \geq 0, \; \phi \in \H,
$$
where $\omega\in\erre$. Then $N$ is $\omega$-$m$-accretive in $\H$,
i.e. 
$$
\ip{N\phi}{\phi}_\H \geq -\omega |\phi|^2_\H
\qquad
\forall \varphi \in D(N)
$$
and $R(\lambda I+N)=\H$ for all $\lambda>\omega$, where $D(\cdot)$ and
$R(\cdot)$ denote domain and range, respectively. Let $g\in\H$ be a
given function and define the closed convex subset of $\H$
$$
\K_g = \{\phi \in \H:\; \phi \geq g \;\mu\textrm{-a.e.} \}.
$$
The normal cone to $\K_g$ at $\phi$ is defined by
$$
\mathcal{N}_g(\phi) = 
\Big\{z \in \H:\, \int_H z(\phi-\psi)\,d\mu \geq 0 \;\;
\forall \psi \in \K \Big\},
$$
or equivalently
$$
\mathcal{N}_g(\phi) = 
\Big\{ z \in \H: \; z(x)=0 \; \textrm{if $\phi(x)>g(x)$},\;
z(x)\leq 0 \; \textrm{if $\phi(x)=g(x)$},\; \mu\textrm{-a.e.}
\Big\}.
$$
We are going to study the parabolic variational inequality
\begin{equation}
  \label{eq:pvi}
\left\{
\begin{array}{ll}
\ds \frac{d\varphi}{dt}(t) + N\varphi(t) + \mathcal{N}_g(\varphi(t)) 
\ni f(t),  & t\in (0,T) \\[8pt]
\varphi(0)=\varphi_0,
\end{array}\right.
\end{equation}
where $\varphi_0 \in \H$ and $f \in L^2([0,T],\H)$ are given.

By a strong solution of (\ref{eq:pvi}) we mean an absolutely
continuous function $\varphi:[0,T]\to\H$ which satisfies
(\ref{eq:pvi}) a.e. on $(0,T)$. A function $\varphi\in
C([0,T],\H)$ is said to be a generalized solution of (\ref{eq:pvi}) if
there exist sequences $\{\varphi_0^n\}\subset\H$, $\{f_n\}\subset
L^2([0,T],\H)$ and $\{\varphi_n\}\subset C([0,T],\H)$ such that, for
all $n$, $\varphi_n$ is a strong solution of
$$
\frac{d\varphi}{dt}(t) + N\varphi(t) + \mathcal{N}_g(\varphi(t)) 
\ni f_n(t)
$$
a.e. on $(0,T)$ with initial condition $\varphi(0)=\varphi_0^n$, and
$\varphi_n\to\varphi$ in $C([0,T],\H)$ as $n\to\infty$.

In order to establish existence of a solution for equation
(\ref{eq:pvi}) we are going to apply the general theory of existence
for Cauchy problems in Hilbert spaces associated with nonlinear maximal
monotone operators (see e.g. \cite{barbu-nonlin}, \cite{barbu},
\cite{Bmax}). We recall that the nonlinear (multivalued) operator
$A:D(A)\subset\H\to\H$ is said to be maximal monotone (or equivalently
$m$-accretive) if $\ip{y_1-y_2}{x_1-x_2} \geq 0$ for all $y_i\in
Ax_i$, $i=1,2$, and $R(I+A)=\H$. The operator $A$ is said to be
$\omega$-$m$-accretive if $\lambda I + A$ is $m$-accretive for all
$\lambda>\omega$. If $A$ is $\omega$-$m$-accretive we set (Yosida
approximation)
$$
A_\lambda u = \frac1\lambda(u-(I+\lambda A)^{-1}u),
\qquad
u\in\H, \;\; 0<\lambda<\frac1\omega.
$$
Recall that $A_\lambda$ is Lipschitz and ${\omega \over
  1-\lambda\omega}$-accretive on $\H$, i.e.
$$
\ip{A_\lambda u - A_\lambda v}{u-v}_\H \geq
- {\omega \over 1-\lambda\omega} |u-v|_\H^2.
$$
Moreover, recalling that $N$ is $\omega$-$m$-accretive, we have the
following result.
\begin{thm}\label{thm:abs}
  Assume that $P_t$ is positivity preserving (that is $P_t\varphi\geq
  0$ for all $\varphi\geq 0$ $\mu$-a.e.) and
\begin{equation}
\label{eq:hyp1}
|(N_\lambda g)^+|_\H \leq C \qquad \forall \lambda\in (0,1/\omega).
\end{equation}
Then the operator $N+\mathcal{N}_g$ with the domain $D(N)\cap\K_g$ is
$\omega$-$m$-accretive in $\H$.
\end{thm}
\begin{proof}
  It is easily seen that $N+\mathcal{N}_g+\omega I$ is accretive. In
  order to prove $m$-accretivity, let us fix $f\in\H$ and
  consider the equation
\begin{equation}
  \label{eq:Gamma}
\alpha\varphi_\lambda + N_\lambda\varphi_\lambda + \mathcal{N}_g(\varphi_\lambda)
\ni f,
\end{equation}
which admits a unique solution for $\alpha>\omega/(1-\lambda\omega)$,
because the operator $N_\lambda+\mathcal{N}_g+\alpha I$ is $m$-accretive for
$\alpha>\omega/(1-\lambda\omega)$.  We are going to show that, as
$\lambda\to 0$, $\varphi_\lambda\to\varphi$ strongly in $\H$ to a
solution $\varphi$ of
\begin{equation}
  \label{eq:ell}
\alpha\varphi + N\varphi + \mathcal{N}_g(\varphi) \ni f.  
\end{equation}
Let us rewrite (\ref{eq:Gamma}) as
\begin{equation}
\label{eq:appr2}
\alpha\psi_\lambda + N_\lambda\psi_\lambda 
+ \mathcal{N}_\K(\psi_\lambda) \ni f - \alpha g -N_\lambda g,
\end{equation}
where $\psi_\lambda = \varphi_\lambda-g$, $\K=\{\psi\in\H:\;\psi\geq 0
\;\; \mu\textrm{-a.e.} \}$, and $\mathcal{N}_\K$ is the normal cone to
$\K$. Setting $\eta_\lambda \in \mathcal{N}_{\K}(\psi_\lambda)$ and multiplying
both sides of (\ref{eq:appr2}) by $\eta_\lambda$ we have
\begin{equation}
\label{eq:appr3}
\alpha\ip{\psi_\lambda}{\eta_\lambda}_\H + |\eta_\lambda|^2_\H
+ \ip{N_\lambda\psi_\lambda+N_\lambda g}{\eta_\lambda}_\H =
\ip{f- \alpha g}{\eta_\lambda}_\H.
\end{equation}
Since $\ip{\psi_\lambda}{\eta_\lambda}_\H \geq 0$ (by definition of
$\mathcal{N}_\K$) and $\ip{N_\lambda\psi_\lambda}{\eta_\lambda}_\H\geq
0$ (in fact $(I+\lambda N)^{-1}\K \subset \K$ because $P_t$ is
positivity preserving), (\ref{eq:appr3}) yields
\begin{equation}
\label{eq:appr4}
|\eta_\lambda|^2_\H + \ip{N_\lambda g}{\eta_\lambda}_\H \leq
\ip{f - \alpha g}{\eta_\lambda}_\H.
\end{equation}
On the other hand, we have $\ip{N_\lambda g}{\eta_\lambda}_\H \geq
\ip{(N_\lambda g)^+}{\eta_\lambda}_\H$, because $\eta_\lambda \in
\mathcal{N}_{\K}(\psi_\lambda)$ implies that
$\ip{\eta_\lambda}{\phi}_\H \leq 0$ if $\phi\geq 0$ $\mu$-a.e..
Then by (\ref{eq:appr4}) and assumption (\ref{eq:hyp1}) we obtain
$$
|\eta_\lambda|_\H \leq |f-\alpha g|_\H + |(N_\lambda g)^+|_\H \leq C
\qquad \forall \lambda \in (0,\omega^{-1}).
$$
Moreover, (\ref{eq:appr2}) implies that
$$
|\psi_\lambda|_\H \leq |f-\alpha g|_\H 
\qquad \forall \lambda \in (0,\omega^{-1}).
$$
Therefore $\{\varphi_\lambda=\psi_\lambda+g\}$ and $\{\eta_\lambda\}$
are bounded in $\H$, and so is $\{N_\lambda\varphi_\lambda\}$.  This
implies by standard arguments that $\{\varphi_\lambda\}$ is Cauchy in
$\H$, so we have that on a subsequence, again denoted by $\lambda$,
$$
\begin{array}{llll}
\varphi_\lambda &\to& \varphi & \textrm{strongly in $\H$},\\
N_\lambda(\varphi_\lambda) &\to& \xi & \textrm{weakly in $\H$},\\
\eta_\lambda &\to& \eta       & \textrm{weakly in $\H$},
\end{array}
$$
as $\lambda \to 0$.  Since
$\eta_\lambda\in\mathcal{N}_g(\varphi_\lambda)$ and $\mathcal{N}_g$ is
maximal monotone, we have $\eta\in\mathcal{N}_g(\varphi)$ and,
similarly, $\xi=N\varphi$. Hence $\varphi$ is a solution of
(\ref{eq:ell}), as required.
\end{proof}
\begin{rmk}
  If $P_t$ is the transition semigroup associated to a Markov
  stochastic process $X$, then $P_t$ is automatically positivity
  preserving. Assumption (\ref{eq:hyp1}) holds in particular if $g\in
  D(N)$ or $(I+\lambda N)^{-1}g \geq g$ for all $\lambda \in
  (0,1/\omega)$.
\end{rmk}
\begin{rmk}
  Denoting by $N^*$ the dual of $N$, the operator $N$ has a natural
  extension from $\H$ to $(D(N^*))'$ defined by
  $Nu(\varphi)=u(N^*\varphi)$ for all $\varphi\in D(N^*)$ and
  $u\in\H$.  Then as $\lambda\to 0$ one has $N_\lambda g \to Ng$
  weakly in $(D(N^*))'$ and if it happens that $Ng$ belongs to a
  lattice subspace, then condition (\ref{eq:hyp1}) simply means that
  $(Ng)^+\in\H$. This is the case in spaces $L^2(\Xi)$,
  $\Xi\subset\erre^n$, where usually $Ng$ is a measure on $\Xi$ (see
  e.g. \cite{Buni}).
\end{rmk}

\begin{rmk}
  Theorem \ref{thm:abs} remains true if we replace assumption
  (\ref{eq:hyp1}) by
  \begin{equation}
    \label{eq:hyp1-alt}
    \frac1t |(g-P_tg)^+|_\H \leq C \qquad \forall t\in(0,1).
  \end{equation}
  The proof follows along completely similar lines.
\end{rmk}

By the general theory of Cauchy problems associated with nonlinear
$m$-accretive operators (see e.g.  \cite{barbu-nonlin}, \cite{barbu},
\cite{Bmax}) we obtain the following result.
\begin{thm}\label{cor:abs}
  Assume that the hypotheses of Theorem \ref{thm:abs} are satisfied. Let
  $\varphi_0 \in D(N)\cap\K_g$ and $f\in W^{1,1}([0,T];\H)$. Then there exists
  a unique strong solution $\varphi \in W^{1,\infty}([0,T];\H) \cap
  L^\infty([0,T];D(N))$ of the Cauchy problem (\ref{eq:pvi}).
  Moreover the function $t \mapsto \varphi(t)$ is right-differentiable and
  $$
  \frac{d^+}{dt}\varphi(t)+\phi(t) = 0, \qquad t\in [0,T),
  $$
  where
  $$
  \phi(t) = \left\{\begin{array}{ll}
  N\varphi(t)-f(t)& \mu\mathrm{-a.e.\ in\ } \{\varphi(t,x)>g(x)\} \\[6pt]
  (N\varphi(t)-f(t))^+& \mu\mathrm{-a.e.\ in\ } \{\varphi(t,x)=g(x)\}.
  \end{array}\right.
  $$
  If $\varphi_0 \in \K_g$ and $f\in L^2([0,T];\H)$ then equation
  (\ref{eq:pvi}) has a unique generalized solution $\varphi\in
  C([0,T],\H)$, $\varphi(t)\in\K_g$ for almost all $t\in [0,T]$.
\end{thm}
We shall see later (see Theorem \ref{prop:gt} below) that the
generalized solution satisfies (\ref{eq:pvi}) in a more precise sense.
\begin{rmk}
  By the general theory of Cauchy problems for nonlinear accretive
  operators (see \cite{barbu-nonlin}, \cite{barbu}, \cite{Bmax}) one
  knows that the solution $\varphi(t)$ given by Theorem \ref{cor:abs}
  can be approximated as $h\to 0$ by the solution
  $\{\varphi_i\}_{i=1}^{N_h}$ of the finite difference scheme
  $$
  \varphi_{i+1} + hN\varphi_{i+1} + h \mathcal{N}_g(\varphi_{i+1})
  \ni f_i + \varphi_i,
  \qquad
  i=0,1,\ldots,N_h,
  $$
  where $hN_h=T$ and $f_i=\int_{ih}^{(i+1)h} f(t)\,dt$. Equivalently,
  $$
  \left\{\begin{array}{ll}
  \varphi_{i+1} = (I+hN)^{-1}(f_i + \varphi_i), &
          \textrm{if } (I+hN)^{-1}(f_i + \varphi_i) > g,\\[4pt]
  \varphi_i > g & \forall i.
  \end{array}\right.
  $$
\end{rmk}

\subsection{Time-dependent obstacle}
We shall consider the case where the obstacle function $g$ depends
also on time. In particular, we shall assume that
\begin{eqnarray}
  \label{eq:gt1}
&&  g \in W^{1,\infty}([0,T],\H)\\
  \label{eq:gt2}
&&  \int_0^T |(N_\lambda g)^+|_\H^2\,dt \leq C
\qquad \forall \lambda \in (0,\omega^{-1}).
\end{eqnarray}
Let $g_\lambda=(I+\lambda N)^{-1}g$ and consider the approximating equation
\begin{equation}
  \label{eq:nurb}
{d\varphi_\lambda \over dt}(t) 
+ N(\varphi_\lambda(t) + g_\lambda(t)-g(t))
+ \mathcal{N}_{g(t)}(\varphi_\lambda(t)) \ni f(t)
\end{equation}
on $(0,T)$ with initial condition $\varphi(0)=\varphi_0$, and
$\varphi_0\geq g(0)$, $f\in L^2([0,T],\H)$. Equivalently, setting
$\psi_\lambda=\varphi_\lambda-g$, we get
\begin{equation}
\label{eq:gt3}
\left\{\begin{array}{l}
\ds {d\psi_\lambda \over dt}(t) + N\psi_\lambda(t)
+ \mathcal{N}_{\K}(\psi_\lambda(t)) 
\ni f(t) - \frac{dg}{dt}(t) - N_\lambda g(t),\\[6pt]
\ds \psi_\lambda(0)=\varphi_0-g(0)\in \K.
\end{array}\right.
\end{equation}
In order to work with strong solutions of equation (\ref{eq:gt3}), we
shall assume, without any loss of generality, that $f\in
W^{1,1}([0,T],\H)$, $\frac{dg}{dt} \in W^{1,1}([0,T],\H)$, and
$\varphi_0-g(0) \in \K \cap D(N)$. This can be achieved in the
argument which follows by taking smooth approximations of $f$, $g$ and
$\varphi_0$.
Then equation (\ref{eq:gt3}) has a unique strong solution
$\psi_\lambda\in W^{1,\infty}([0,T],\H)\cap L^\infty([0,T],D(N))$ by
standard existence results for Cauchy problems because, as seen
earlier, $N+\mathcal{N}_\K$ is $\omega$-$m$-accretive.  Moreover,
multiplying both sides of (\ref{eq:gt3}) by
$\eta_\lambda(t)\in\mathcal{N}_\K(\psi_\lambda(t))$ and taking into
account that $P_t$ is positivity preserving and 
$$
\int_0^T \ip{N\psi_\lambda}{\eta_\lambda}_\H dt \geq 0,
\qquad
\int_0^T \ip{\frac{d\psi_\lambda}{dt}(t)}{\eta_\lambda(t)}_\H\,dt = 0
\qquad
\forall \lambda \in (0,\omega^{-1}),
$$
arguing as in the proof of Theorem \ref{thm:abs}, we get the
following a priori estimates:
\begin{eqnarray}
  \label{eq:monza}
  |\varphi_\lambda(t)|_\H &\leq& C \qquad \forall t\in[0,T], \\
  \label{eq:sm}
  \int_0^T|\eta_\lambda(t)|^2_\H\,dt &\leq& C,
\end{eqnarray}
for all $\lambda\in(0,\omega^{-1})$.
Hence on a subsequence, again denoted by $\lambda$, we have
$$
\begin{array}{llll}
\varphi_\lambda &\to& \varphi & \textrm{weakly* in $L^\infty([0,T],\H)$} \\
\eta_\lambda &\to& \eta& \textrm{weakly in $L^2([0,T],\H)$}
\end{array}
$$
as $\lambda \to 0$.
Moreover, $\varphi:[0,T]\to\H$ is weakly continuous and
\begin{equation}
  \label{eq:wc}
\frac{d\varphi}{dt}(t)+N\varphi(t)+\eta(t) = f(t)  
\end{equation}
almost everywhere in $[0,T]$ with initial condition
$\varphi(0)=\varphi_0$ in mild sense, i.e.,
$$
\varphi(t) + \int_0^t e^{-N(t-s)}\eta(s)\,ds =
e^{-Nt}\varphi_0 + \int_0^t e^{-N(t-s)} f(s)\,ds
$$
for almost all $t\in [0,T]$. The latter follows by letting $\lambda \to 0$ into
the equation
\begin{equation}
\label{eq:16}
\begin{array}{c}
\ds \varphi_\lambda(t) + g_\lambda(t) - g(t) 
+ \int_0^t e^{-N(t-s)}(\eta_\lambda(s)-f(s)-g'_\lambda(s)+g'(s))\,ds\\[8pt]
\ds \qquad\qquad\qquad\qquad = e^{-Nt}(\varphi_0 + g_\lambda(0) - g(0)).
\end{array}
\end{equation}
Taking into account that, as $\lambda \to 0$, $g_\lambda(t) \to g(t)$
strongly in $\H$ on $[0,T]$ and $g'_\lambda-g'=(I+\lambda N)^{-1}g'-g'
\to 0$ strongly in $L^2([0,T],\H)$, we obtain the desired equation.
In particular it follows that $\varphi_\lambda(t) \to \varphi(t)$
weakly in $\H$ for $t\in[0,T]$.
We are going to show that $\eta(t) \in \mathcal{N}_g(\varphi(t))$ a.e.
on $[0,T]$. To this purpose it suffices to show that
\begin{equation}
\label{eq:lino}
\limsup_{\lambda\to 0} \int_0^T
e^{\gamma t} \ip{\eta_\lambda(t)}{\varphi_\lambda(t)}_\H\,dt
\leq
\int_0^T e^{\gamma t} \ip{\eta(t)}{\varphi(t)}_\H\,dt
\end{equation}
for some real number $\gamma$. We shall prove that (\ref{eq:lino})
holds with $\gamma=-2\omega$.
To this end we set $N_\omega=N+\omega I$ (note that $N_\omega$ is
$m$-accretive in $\H$) and, rewriting equation (\ref{eq:nurb}) as
$$
\frac{d}{dt}(\varphi_\lambda + g_\lambda - g) 
+ N_\omega(\varphi_\lambda + g_\lambda - g) 
+ \eta_\lambda - \omega(\varphi_\lambda + g_\lambda - g) 
= f + g'_\lambda - g',
$$
we may equivalently write (\ref{eq:16}) as
$$
\begin{array}{c}
\ds e^{-\omega t} (\varphi_\lambda(t) + g_\lambda(t) - g(t))
+ \int_0^t e^{-N_\omega(t-s)} e^{-\omega s}
      (\eta_\lambda(s)-f(s)-g'_\lambda(s)+g'(s))\,ds\\[8pt]
\ds \qquad\qquad\qquad\qquad = e^{-N_\omega t}(\varphi_0 + g_\lambda(0) - g(0))
\qquad \forall t \in (0,T).
\end{array}
$$
This yields
\begin{eqnarray*}
&&\int_0^T e^{-2\omega t}\ip{\eta_\lambda(t)}{\varphi_\lambda(t)}_\H\,dt =\\
&&\qquad -\int_0^T\ip{e^{-\omega t}\eta_\lambda(t)}
            {\int_0^t e^{-N_\omega(t-s)} e^{-\omega s}\eta_\lambda(s)}_\H ds\\
&&\qquad + \int_0^T e^{-\omega t}\ip{\eta_\lambda(t)}
   {e^{-N_\omega t}(\varphi_0+g_\lambda(0)-g(0))
       -e^{-\omega t}(g_\lambda(t)-g(t))}_\H dt\\
&&\qquad + \int_0^T e^{-\omega t}\ip{\eta_\lambda(t)}
   {\int_0^t e^{-N_\omega(t-s)} e^{-\omega s} (f(s)+g'_\lambda(s)-g'(s))\,ds}_\H dt
\end{eqnarray*}
Then letting $\lambda\to 0$ we obtain
\begin{equation}
\label{eq:hock}
\begin{array}{l}
\ds \limsup_{\lambda\to 0}
\int_0^T e^{-2\omega t} \ip{\eta_\lambda(t)}{\varphi_\lambda(t)}_\H\,dt \leq\\[8pt]
\qquad\ds -\liminf_{\lambda\to 0}
   \int_0^T\ip{e^{-\omega t}\eta_\lambda(t)}
      {\int_0^t e^{-N_\omega(t-s)} e^{-\omega s}\eta_\lambda(s)\,ds}_\H dt\\[8pt]
\qquad\ds + \int_0^T e^{-\omega t}\ip{\eta(t)}
    {e^{-N_\omega t}\varphi_0+\int_0^te^{-N_\omega(t-s)} e^{-\omega s}f(s)\,ds}_\H dt.
\end{array}
\end{equation}
Consider the function
$$
F(y) = \int_0^T \ip{y(t)}{\int_0^t e^{-N_\omega(t-s)}y(s)\,ds}_\H\,dt,
\qquad y \in L^2([0,T],\H),
$$
which is continuous and convex on $L^2([0,T],\H)$ (the latter is an
easy consequence of the fact that $N_\omega$ is accretive). Hence $F$
is weakly lower semicontinuous and therefore
$$
\liminf_{\lambda\to 0} F(e^{-\omega t}\eta_\lambda) \geq F(e^{-\omega t}\eta).
$$
Substituting this expression into (\ref{eq:hock}) we find that
\begin{equation}
\label{eq:19}
\begin{array}{l}
\ds \limsup_{\lambda\to 0}
\int_0^T e^{-2\omega t} \ip{\eta_\lambda(t)}{\varphi_\lambda(t)}_\H\,dt \leq\\[8pt]
\qquad\ds - \int_0^T\ip{e^{-\omega t}\eta(t)}
      {\int_0^t e^{-N_\omega(t-s)} e^{-\omega s}\eta(s)\,ds}_\H dt\\[8pt]
\qquad\ds + \int_0^T e^{-\omega t}\ip{\eta(t)}
    {e^{-N_\omega t}\varphi_0
         +\int_0^te^{-N_\omega(t-s)} e^{-\omega s}f(s)\,ds}_\H dt\\[8pt]
\qquad\ds = \int_0^T e^{-2\omega t}\ip{\eta(t)}{\varphi(t)}_\H dt.
\end{array}
\end{equation}
The latter follows by equation $d\varphi/dt+N\varphi+\eta=f$, or equivalently
$$
\frac{d}{dt}(e^{-\omega t}\varphi(t)) + N_\omega(e^{-\omega t}\varphi(t))
+\eta(t)e^{-\omega t} = e^{-\omega t} f(t).
$$
Hence $\eta(t) \in \mathcal{N}_g(\varphi(t))$ for all $t\in (0,T)$ as
claimed.

\begin{defi}
A function $\varphi\in C([0,T],\H)$ is said to be a mild solution of 
\begin{equation}
  \label{eq:nurba}
{d\varphi \over dt}(t) + N\varphi(t)
+ \mathcal{N}_{g(t)}(\varphi(t)) \ni f(t)
\end{equation}
on $[0,T]$ with initial condition $\varphi(0)=\varphi_0$ if
$\varphi(t)\geq g(t)$ $\mu$-a.e. for almost all $t\in[0,T]$ and
there exists $\eta\in L^2([0,T],\H)$ with
$\eta(t)\in\mathcal{N}_{g(t)}(\varphi(t))$ for almost all $t\in[0,T]$, such
that
\begin{equation}
  \label{eq:mild}
\varphi(t) + \int_0^t e^{-N(t-s)}\eta(s)\,ds =
e^{-Nt}\varphi_0 + \int_0^t e^{-N(t-s)} f(s)\,ds
\end{equation}
for all $t\in[0,T]$.
\end{defi}
\begin{thm}\label{prop:gt}
  Assume that $P_t$ is positivity preserving and (\ref{eq:gt1}),
  (\ref{eq:gt2}) hold. Let $\varphi_0\in\H$, $\varphi_0\geq g(0)$ and
  $f\in L^2([0,T],\H)$. Then (\ref{eq:nurba}) has a unique mild
  solution. Moreover, the map $(\varphi_0,f) \mapsto \varphi$ is
  Lipschitz from $\H\times L^2([0,T],\H)$ to $C([0,T],\H)$.
\end{thm}
\begin{proof}
  Existence was proved above. Uniqueness as well as as continuous
  dependence on data follows by (\ref{eq:mild}) taking into account
  that $\eta(t)\in\mathcal{N}_{g(t)}(\varphi(t))$ for almost all $t\in[0,T]$
  and
\begin{eqnarray*}
\int_0^T \ip{\eta(t)}{\int_0^t e^{-N(t-s)}\eta(s)\,ds}_\H dt 
&\geq& -\omega \int_0^T \left|\int_0^t e^{-N(t-s)}\eta(s)\,ds\right|_\H^2 dt\\
&& + \frac12 \left|\int_0^T e^{-N(t-s)}\eta(s)\,ds \right|_\H^2.
\end{eqnarray*}
\end{proof}
It is worth emphasizing that in the case where $g$ is time-dependent
the ``mild'' solution provided by Theorem \ref{prop:gt} is a
generalized solution in the sense of Theorem \ref{cor:abs}. However,
even in this case Theorem \ref{prop:gt} is not directly implied by
Theorem \ref{cor:abs}.

\section{Variational inequalities and optimal stopping
problems}
\label{sec:viosp}
\setcounter{equation}{0} Let $H$ be a real Hilbert space with inner
product $\ip{\cdot}{\cdot}$ and norm $|\cdot|$, and
$(\Omega,\mathcal{F},\mathbb{F}=(\mathcal{F}_t)_{t\geq 0},\mathbb{P})$
a filtered probability space satisfying the usual conditions, on which
an $H$-valued Wiener process (adapted to $\mathbb{F}$) with covariance
operator $Q$ is defined. Let $X$ be the process generated by the
stochastic differential equation
\begin{equation}
\label{eq:pierogi}
dX(s) = b(X(s))\,ds + \sigma(X(s))\,dW(s)  
\end{equation}
on $s\in[t,T]$ with initial condition $X(t)=x$, where $b:H\to H$ and
$\sigma:H\to L(H,H)$ are such that (\ref{eq:pierogi}) admits a unique
solution that is strong Markov.  Define the value function $v(t,x)$ of
an optimal stopping problem for $X$ as
\begin{equation}
\label{eq:fval}
v(t,x) = \sup_{\tau\in\mathfrak{M}} \E_{t,x}\Big[
e^{-\psi(t,\tau)} g(\tau,X(\tau))
+ \int_t^\tau e^{-\psi(t,s)}f(s,X(s))\,ds
\Big],
\end{equation}
where $\mathfrak{M}$ is the family of all $\mathbb{F}$-stopping times
such that $\tau\in[t,T]$ $\mathbb{P}$-a.s., and
$$
\psi(t,s) = \int_t^s c(X_r)\,dr \qquad \forall t \leq s \leq T,
$$
where $c:H\to\erre_+$ is a given discount function (which we also
assume to be bounded, for simplicity). Exact conditions on $g$ and $f$
will be specified below. The function $v$ is formally the solution of
the backward variational inequality
\begin{equation}
  \label{eq:23}
  \frac{\partial u}{\partial t} + L_0 u - c u
  - \mathcal{N}_{g(t)}(u) \ni f
\end{equation}
in $(0,T)\times H$ with terminal condition $u(T,x)=g(T,x)$, where
\begin{equation}
  \label{eq:no}
  L_0\phi = \frac12 \tr [(\sigma Q^{1/2})(\sigma Q^{1/2})^* D^2\phi]
         + \ip{b(x)}{D\phi}, \qquad \phi \in D(L_0) = C^2_b(H).
\end{equation}
More precisely, denoting by $\mu$ an excessive measure of the
transition semigroup $P_t$ generated by the process $X$,
we have that $v$ is the solution of the variational inequality
\begin{equation}
  \label{eq:roga}
  \frac{\partial u}{\partial t} + L u - cu
  - \mathcal{N}_{g(t)}(u) \ni f
\end{equation}
in $(0,T)$ with terminal condition $u(T)=g(T)$, where $L$ is the
infinitesimal generator of $P_t$. In many situations of interest
$L=\overline{L_0}$, the closure of $L_0$ in $L^2(H,\mu)$.  Before
giving a simple sufficient condition for this to hold, let us
define precisely excessive measures.
\begin{defi}
  Let $P_t$ be a strongly continuous semigroup on $L^2(H,\mu)$, where
  $\mu$ is a probability measure on $H$. The measure $\mu$ is called
  excessive for $P_t$ if there exists $\omega>0$ such that
  \[
    \int_H P_tf \,d\mu \leq e^{\omega t} \int_H f\,d\mu
    \qquad \forall t\geq 0
  \]
  for all bounded Borel functions $f$ with $f \geq 0$ $\mu$-a.e..
\end{defi}


We have then the following result.
\begin{lemma}\label{lem:clos}
  Let $P_t$ the semigroup generated by $X$, and let $\mu$ be an
  excessive measure for $P_t$ on $H$. Moreover, let $b \in
  C^2(H) \cap L^2(H,\mu)$, $\sigma\in C^2(H,L(H,H))$, and
  \begin{equation}
    \label{eq:bsig}
    |Db(x)|_H + |D\sigma(x)|_{L(H,H)} \leq C
  \end{equation}
  for all $x\in H$. Then $-L_0$ is $\omega$-accretive and $L$ is the
  closure in $L^2(H,\mu)$ of $L_0$ defined on $D(L_0)=C^2_b(H)$.
\end{lemma}
\begin{proof}
  The argument is similar to that used in \cite{DP-K} for similar problems, so
  it will be sketched only. Fix $h \in C^2_b(H)$ and consider the
  equation $(\lambda I - L_0)\varphi=h$, or equivalently
  \begin{equation}
    \label{eq:resolv}
  \varphi(x) = \E_{0,x}\int_0^\infty e^{-\lambda t} h(X(t))\,dt,
  \qquad \lambda > \omega.
  \end{equation}
It is readily seen that $\varphi\in C^2_b(H)$ and, by It\^o's formula,
$(\lambda-L_0)\varphi=h$ in $H$. Since $-L_0$ is closable and
$\omega$-accretive, and $R(\lambda-L_0)$ is dense in $L^2(H,\mu)$, we
infer that
$\overline{L_0}$ coincides with $L$.
\end{proof}
Note also that since the measure $\mu$ is $\omega$-excessive for $P_t$
we have $\int_H Lf\,d\mu \leq \omega \int_H f\,d\mu$, which implies
that $L$ is $\omega$-dissipative in $L^2(H,\mu)$. In the sequel, for
convenience of notation, we shall set $N=-L+cI$.

We shall further assume that $g(t,x)$ is continuously differentiable with
respect to $t$, Lipschitz in $x$, and
\begin{eqnarray}
&& \sup_{t\in(0,T)} \int_H (|D_tg(t,x)|^2+|D_xg(t,x)|^2)\,\mu(dx) 
     < \infty,\label{eq:poco}\\[4pt]
&& \tr[(\sigma Q^{1/2})(\sigma Q^{1/2})^* D_{xx}^2g] \geq 0
\qquad
\textrm{on\ } (0,T)\times H.
\label{eq:troppo}
\end{eqnarray}
If $H$ is a finite dimensional space, the inequality (\ref{eq:troppo})
must be interpreted in the sense of distributions (i.e. of measures).
In the general situation treated here the exact meaning of
(\ref{eq:troppo}) is the following: there exists a sequence
$\{g_\varepsilon(t)\}\subset C_b^2(H)$ such that
$$
\begin{array}{ll}
\ds \sup_{t\in(0,T)} \int_H (|D_tg_\varepsilon(t,x)|^2
           +|D_xg_\varepsilon(t,x)|^2)\,\mu(dx) < C
     & \forall \varepsilon>0,\\[14pt]
\ds \tr[(\sigma Q^{1/2})(\sigma Q^{1/2})^* D_{xx}^2g_\varepsilon(t,x)]
       \geq 0 & \forall \varepsilon>0,\; t\geq0,\;x\in H,\\[10pt]
\ds g_\varepsilon(t) \to g(t) & \textrm{in\ } L^2(H,\mu)
                                      \;\; \forall t\geq 0.
\end{array}
$$
It turns out that under assumption (\ref{eq:troppo}) $g$ satisfies
condition (\ref{eq:gt2}). Here is the argument: for each $\lambda>0$
we have $(N_\lambda g)^+=\lim_{\varepsilon\to 0} (N_\lambda
g_\varepsilon)^+$ in $L^2(H,\mu)$. On the other hand, $N_\lambda
g_\varepsilon = N(I+\lambda N)^{-1}g_\varepsilon$ and by
(\ref{eq:troppo}) we see that
$$
\tr\Big[
    (\sigma Q^{1/2})(\sigma Q^{1/2})^*
    D_{xx}^2[(I+\lambda N)^{-1} g_\varepsilon]
\Big] \geq 0 \qquad \textrm{on\ } H
$$
because $(I+\lambda N)^{-1}$ leaves invariant the cone of nonnegative
functions (by the positivity preserving property of $P_t$). Hence
$$
\Big|(N_\lambda g_\varepsilon)^+\Big|_{L^2(H,\mu)} \leq
\Big| \ip{b}{D_x(I+\lambda N)^{-1}g_\varepsilon} \Big|_{L^2(H,\mu)} 
\leq C \quad\;\; \forall \lambda\in(0,\omega^{-1}),\,\varepsilon>0
$$
because $b\in L^2(H,\mu)$. This implies (\ref{eq:gt2}) as claimed.
\begin{prop}\label{prop:brutta}
  Assume that $f\in L^2([0,T],L^2(H,\mu) \cap C([0,T],C_b(H))$ and
  that conditions (\ref{eq:gt1}), (\ref{eq:bsig}) and
  (\ref{eq:troppo}) hold. Furthermore, assume that the law of $X(s)$
  is absolutely continuous with respect to $\mu$ for all $s\in[t,T]$.
  Then there exists a unique mild solution $u\in C([0,T];L^2(H,\mu))$
  of the variational inequality (\ref{eq:roga}).  Moreover, $u$
  coincides $\mu$-a.e. with the value function $v$ defined in
  (\ref{eq:fval}).
\end{prop}
\begin{proof}
  Existence and uniqueness for (\ref{eq:roga}) follows by Proposition
  \ref{prop:gt}. In the remaining of the proof we shall limit
  ourselves to the case $f=0$. This is done only for simplicity, as
  the reasoning is identical in the more general case $f\neq 0$.
  By definition of mild solution there exists $\eta\in
  L^2([0,T],L^2(H,\mu))$ such that
  $\eta(t)\in\mathcal{N}_{g(t)}(u(t))$ for all $t\in[0,T]$ and the
  following equation is satisfied (in mild sense) for all $s\in(0,T)$,
  with terminal condition $u(T)=g(T)$:
\begin{equation}
\label{eq:mimi}
\frac{du}{ds}(s) - Nu(s) = \eta(s),
\end{equation}
i.e.,
\begin{equation}
  \label{eq:metal}
u(t,x) = -\int_t^\tau R_{s-t}\eta(s,x)\,ds + R_{\tau-t}u(\tau,x)
\quad
\forall t < \tau < T, \; \mu\textrm{-a.e.}\; x \in H,
\end{equation}
where $R_t$ is the transition semigroup generated by $-N$, or
equivalently the following Feynman-Kac semigroup associated with the
stochastic differential equation (\ref{eq:pierogi}):
$$
R_t\phi(x) = \E_{0,x}\Big[ e^{-\int_0^t c(X(s))\,ds} \phi(X(t)) \Big],
\qquad \phi \in L^2(H,\mu).
$$
Let us set $H_T=[t,T]\times H$ and define the measure
$\mu_T=\mathrm{Leb}\times\mu$ on $H_T$, where $\mathrm{Leb}$ stands
for one-dimensional Lebesgue measure. Recalling that $u(s,x) \geq
g(s,x)$ for all $s\in [t,T]$, $\mu$-a.e. $x\in H$, we can obtain a
version of $u$, still denoted by $u$, such that $u(s,x)\geq g(s,x)$
for all $(s,x)\in H_T$. Recalling that $\eta(s,\cdot) \in L^2(H,\mu)$
a.e. $s\in [t,T]$, equation (\ref{eq:metal}) yields
\begin{equation}
  \label{eq:lurgico}
u(t,x) = \E_{t,x} \Big[ \int_t^\tau -e^{-\psi(t,s)} \eta(s,X(s))\,ds
+ e^{-\psi(t,\tau)} u(\tau,X(\tau)) \Big]
\end{equation}
for every stopping time $\tau\in[t,T]$, for all $t\in[0,T]$ and
$\mu$-a.e. $x \in H$. In fact, let us consider a sequence
$\{\eta_\varepsilon\}\subset C^1([0,T],\H)$ such that
$\eta_\varepsilon \to \eta$ in $L^2([0,T],\H)$. Then equation
(\ref{eq:mimi}), with $\eta_\varepsilon$ replacing $\eta$, admits a
solution $u_\varepsilon \in C^1([0,T],\H) \cap C([0,T],D(N))$ such
that $u_\varepsilon \to u$ in $C([0,T],\H)$ as $\varepsilon \to 0$.
Recalling that $N=-\overline{L_0}+cI$, there exists a sequence
$\{w_\varepsilon\} \subset C^1([0,T],\H) \cap C([0,T],C^2_b(H))$ such
that
\begin{eqnarray*}
|u_\varepsilon(t) - w_\varepsilon(t)|_\H &\leq& \varepsilon\\
|Nu_\varepsilon(t) - (-L_0+cI)w_\varepsilon(t)|_\H &\leq& \varepsilon\\
|\frac{du_\varepsilon}{dt}(t) - \frac{dw_\varepsilon}{dt}(t)|_\H &\leq&
\varepsilon
\end{eqnarray*}
for all $t\in[0,T]$. Therefore we have
$$
\frac{dw_\varepsilon}{dt}(t) - (L_0-cI)w_\varepsilon(t) = 
\tilde{\eta}_\varepsilon \qquad \forall t \in [0,T],
$$
where $\tilde{\eta}_\varepsilon \to \eta$ in $L^2([0,T],\H)$. Then we have
\begin{equation}
  \label{eq:lurido}
w_\varepsilon(t,x) = \E_{t,x} \Big[ \int_t^\tau -e^{-\psi(t,s)}
\tilde{\eta}_\varepsilon(s,X(s))\,ds 
+ e^{-\psi(t,\tau)} w_\varepsilon(\tau,X(\tau)) \Big]
\end{equation}
for all stopping times $\tau\in[t,T]$. We shall now show that
(assuming, without loss of generality, $\psi\equiv 0$)
$$
\E_{t,x} \int_t^\tau \tilde{\eta}_\varepsilon(s,X(s))\,ds \to
\E_{t,x} \int_t^\tau \eta(s,X(s))\,ds
$$
for all $t\in[0,T]$ and in $L^2(H,\mu)$ w.r.t. $x$. In fact, Tonelli's
theorem yields, recalling that $\mu$ is excessive for $P_t$,
\begin{eqnarray*}
\lefteqn{\int_H \E_{t,x}\int_t^\tau |\tilde{\eta}_\varepsilon(s,X(s))
   - \eta(s,X(s))|^2\,ds \, \mu(dx)}\\
&\leq& \int_0^T \int_H \E_{0,x} |\tilde{\eta}_\varepsilon(s,X(s))
   - \eta(s,X(s))|^2\,ds \, \mu(dx)\\
&=& \int_0^T \int_H P_s |\tilde{\eta}_\varepsilon(s,x)
   - \eta(s,x)|^2 \, \mu(dx) \, ds \\
&\leq& e^{\omega T} \int_0^T \int_H |\tilde{\eta}_\varepsilon(s,x)
   - \eta(s,x)|^2 \, \mu(dx) \, ds \to 0
\end{eqnarray*}
as $\varepsilon \to 0$, because $\tilde\eta_\varepsilon\to\eta$ in
$L^2([0,T],L^2(H,\mu))$.  An analogous argument shows that
$\E_{t,x}w_\varepsilon(\tau,X(\tau)) \to \E_{t,x}u(\tau,X(\tau))$ for
all $t\in [0,T]$ and in $L^2(H,\mu)$ w.r.t. $x$. Therefore,
passing to a subsequence of $\varepsilon$ if necessary, we have that
the left-hand and right-hand side of (\ref{eq:lurido}) converge
to the left-hand and right-hand side, respectively, of
(\ref{eq:lurgico}) for all $t\in[0,T]$ and $\mu$-a.e. $x\in H$.
Recalling that
$$
\eta(s,x)
\left\{
\begin{array}{lll}
= 0 & \textrm{if\ } u(s,x) > g(s,x) &
\textrm{for each $s$ and $\mu$-a.e. $x\in H$}, \\
\leq 0 & \textrm{if\ } u(s,x) = g(s,x) &
\textrm{for each $s$ and $\mu$-a.e. $x\in H$},
\end{array}\right.
$$
let us define the set
$$
A = \{(s,x) \in H_T: \; \eta(s,x) > 0 \},
$$
for which we have $\mu_T(A)=0$.  Using this fact together with the
assumption that the law of $X(s)$ is absolutely continuous w.r.t.
$\mu$ for all $s\in [t,T]$, hence that $\mathbb{P}_{t,x}((s,X(s))\in
A)=0$, we get
$$
\E_{t,x} \int_t^\tau -e^{-\psi(t,s)} \eta(s,X(s))\,ds \geq 0.
$$
Therefore equation (\ref{eq:lurgico}) implies that $u(t,x)\geq
\E_{t,x}[e^{-\psi(t,\tau)}g(\tau,X(\tau))]$ for all stopping times
$\tau\in\mathfrak{M}$, hence $u(t,x)\geq v(t,x)$, for all $t\in[0,T]$
and $\mu$-a.e. $x\in H$.
Let us now prove that there exists a stopping time $\bar\tau \in
[t,T]$ such that
$u(t,x)=\E_{t,x}[e^{-\psi(t,\bar\tau)}g(\bar\tau,X(\bar\tau))]$, which
will yield $u(t,x)=v(t,x)$, for all $t\in [0,T]$, $\mu$-a.e. $x\in H$.
Define the set
$$
B = \{ (s,x) \in H_T: g(s,x) = u(s,x) \}
$$
and the random time
$$
D_B = \inf\{ s\geq t: (s,X(s)) \in B \} \wedge T.
$$
Since $B$ is a Borel subset of $H_T$ and the process $(s,X(s))$ is
progressive (because it is adapted and continuous), the d\'ebut
theorem (see T.IV.50 in \cite{DM-1}) implies that $D_B$ is a stopping
time. Recalling that $u(s,x)>g(s,x)$ for all $s\in[t,D_B)$, we have,
reasoning as before, $\eta(s,X(s))=0$ a.s. for each $s\in[t,D_B)$.
Thus, taking $\bar\tau=D_B$, (\ref{eq:lurgico}) yields
$$
u(t,x) = \E_{t,x}[e^{-\psi(t,\bar\tau)} g(\bar\tau,X(\bar\tau))]
\qquad \forall t\in[0,T], \; \mu\textrm{-a.e.} \; x\in H.
$$
We have thus proved that there exists a version of $u$ such that
$u(t,x)=v(t,x)$ for all $t\in[t,T]$, $\mu$-a.e. $x\in H$. The
definition of mild solution then implies that $u(t,x)=v(t,x)$ for all
$t\in [t,T]$ and $\mu$-a.e. $x\in H$.
\end{proof}
\begin{rmk}
  The absolute continuity assumption in proposition \ref{prop:brutta}
  can be difficult to verify in general. However, it holds in many
  cases of interest. In particular, it is automatically satisfied if
  the semigroup $P_t$ is irreducible and $\mu$ is invariant with
  respect to $P_t$. Moreover, in the finite dimensional case, if the
  excessive measure $\mu$ is absolutely continuous with respect to
  Lebesgue measure and the coefficients of (\ref{eq:pierogi}) satisfy
  an hypoellipticity condition, the assumptions of the above
  proposition are also satisfied. We shall see in the next section
  that $\mu$ has full support in all examples considered. Moreover, in
  the finite dimensional cases, $\mu$ can be chosen absolutely
  continuous with respect to Lebesgue measure. Let us also remark that
  the continuity of the value function has been proved under very mild
  assumptions by Krylov \cite{krylov}, and by Zabczyk \cite{Zosp} in
  the infinite dimensional case.
\end{rmk}
\begin{rmk}
  Optimal stopping problems in Hilbert spaces and corresponding
  variational inequalities are studied by G{\c{a}}tarek and
  {\'S}wi{\c{e}}ch \cite{GaSw} in the framework of viscosity
  solutions. Their results are applied to pricing interest-rate
  American options, for which the natural dynamics is infinite
  dimensional (e.g. when choosing as state variable the forward
  curve). At the expense of assuming (\ref{eq:troppo}), that is,
  roughly speaking, a convexity assumption on the payoff function $g$,
  we obtained here a more regular solution. We would like to remark
  that $g$ is convex in practically all examples of interest arising
  in option pricing, some of which are investigated in the next section.
\end{rmk}

\section{Pricing of American options}
\setcounter{equation}{0}
Let $\mathbb{Q}$ be a risk neutral martingale measure, and assume we
have $n$ assets whose price-per-share $X(t)=(X_i(t))_{i=1,\ldots,n}$
evolves according to the following Markovian stochastic differential
equation:
\begin{equation}
  \label{eq:sde}
dX(t) = rX(t)\,dt + \sigma(X(t))\,dW(t), \quad X(0)=x\geq 0,
\quad t\in[0,T],
\end{equation}
where $r\in\erre_+$ is the risk-free interest rate, $W$ is a
$\erre^m$-valued Wiener process, and $\sigma:\erre^n\to
L(\erre^m,\erre^n)$ is the volatility function.  Moreover, we assume
that $\sigma$ is such that $X(t)\in\erre^n_+$ for all $t\in[0,T]$. The
standard assumption (see e.g. \cite{K97}) is that
$\sigma_{ij}(X(t))=X_i(t)\tilde\sigma_{ij}(X(t))$ for some
$\tilde\sigma:\erre^n\to L(\erre^m,\erre^n)$. We do \emph{not} assume,
however, that $\sigma$ nor $\tilde\sigma$ satisfies a uniform
nondegeneracy condition. Note that in this situation the market is
incomplete, even if $m=n$, and the choice of the risk neutral
measure $\mathbb{Q}$ is not unique (\cite{K97}).

It is well known that the problem of pricing an American contingent
claim with payoff function $g:\erre^n\to\erre$ is equivalent to the
optimal stopping problem
\begin{equation}  \label{eq:aop}
v(t,x)=\sup_{\tau\in\mathfrak{M}} \E_{t,x}[e^{-r\tau}g(X(\tau))],
\end{equation}
where $\mathfrak{M}$ is the set of all $\mathbb{F}$-adapted stopping
times $\tau\in[t,T]$ and $\E$ stands for expectation with respect to
the measure $\mathbb{Q}$.  Denote by $P_t$ the transition semigroup
associated with (\ref{eq:sde}), i.e. $P_tf(x)=\E_{0,x}f(X(t))$, $f\in
C_b(\erre^n)$, $x\in\erre^n$, and let $L_0$ be the corresponding
Kolmogorov operator. A simple calculation based on It\^o's formula yields
$$
L_0f(x) = \frac12\tr[\sigma(x)\sigma^*(x)D^2f(x)] + \ip{rx}{Df(x)}_{\erre^n},
\qquad
f \in C^2_b(\erre^n).
$$

By classical results (see e.g. \cite{krylov}), the value function
$v(t,x)$ is expected to satisfy the following backward variational
inequality
\begin{equation}
\label{eq:vi}
\left\{
\begin{array}{ll}
\ds \max\Big( (\partial_t+L_0)v(t,x)-rv(t,x),g(x)-v(t,x) \Big)=0,
& (t,x)\in Q_T\\
v(T,x) = g(x), & x\in\erre^n_+,
\end{array}\right.
\end{equation}
where $Q_T=[0,T]\times\erre^n_+$.

The classical theory of variational inequalities in Sobolev spaces
with respect to Lebesgue measure does not apply, however, mainly
because the volatility coefficient is degenerate (see \cite{JLL}).
Nonetheless, one might try to study (\ref{eq:vi}) in spaces of
integrable functions with respect to a suitably chosen measure. The
most natural choice would be an (infinitesimally) invariant measure for
$L_0$. However, without non-degeneracy conditions for $\sigma$ and
with $r>0$, one may not expect existence of an invariant measure (see
e.g. \cite{ABR}, \cite{BKR}).  Here we shall instead solve (\ref{eq:vi}) in
$L^2(\erre^n,\mu)$, where $\mu$ is an (infinitesimally) excessive
measure for $L_0$, which is also absolutely continuous with respect to
Lebesgue measure.

\smallskip

The backward variational inequality (\ref{eq:vi}) can be equivalently
written as the (abstract) variational inequality in
$L^2(\erre^n,\mu)$
\begin{equation}
\label{eq:via}
\partial_t v - Nv - \mathcal{N}_g(v) \ni 0,
\quad v(T)=g,
\end{equation}
where $N=-L+rI$, with $L$ the generator of $P_t$ (which will often
turn out to be the closure of $L_0$), and $\mathcal{N}_g$ is the
normal cone to
$$
\K_g=\{\phi \in L^2(\erre^n,\mu): \; \phi \geq g \;\, \mu\textrm{-a.s} \}.
$$
\begin{lemma}
Assume that
\begin{equation}
\label{eq:hyp-sigma}
\sigma \in C^2(\erre^n), \quad 
|\sigma(x)| \leq C(1+|x|), \quad
|\sigma_{x_i}| + |\sigma_{x_ix_j}| \leq C.
\end{equation}
Then there exists an excessive probability measure $\mu$ of $P_t$ of the form
$$
\mu(dx) = \frac{a}{1+|x|^{2(n+1)}}\,dx
$$
with $a>0$.
\end{lemma}
\begin{proof}
  Setting $\rho(x) = \frac{1}{1+|x|^{2(n+1)}}$, we shall check that
  $L_0^*\rho \leq \omega \rho$ in $\erre^n$ for some $\omega>0$, where
  $L_0^*$ is the formal adjoint of $L_0$, i.e.
  $$
  L_0^*\rho = \frac12\tr[D^2(\sigma\sigma^*\rho)] - r\,\mathrm{div}(x\rho).
  $$
  Assumption (\ref{eq:hyp-sigma}) implies, after some computations, that
  $$
  \sup_{x\in\erre^n} \frac{L_0^*\rho}{\rho} =: \omega < \infty,
  $$
  thus $\mu(dx)=a \rho(x)\,dx$, with $a^{-1}=\int_{\erre^n}
  \rho(x)\,dx$, is a probability measure and satisfies $L_0^*\mu \leq
  \omega\mu$. This yields
  \begin{equation}\label{eq:ex}
  \int_{\erre^n} L_0 f\,d\mu \leq \omega\int_{\erre^n} f \,d\mu
  \end{equation}
  for all $f \in C^2_b(\erre^n)$ with $f\geq 0$, and therefore
  $$
  \int_{\erre^n} P_t f\,d\mu \leq e^{\omega t} \int_{\erre^n} f \,d\mu
  $$
  for all $f \in C^2_b(\erre^n)$, $f\geq 0$. The latter extends by
  continuity to all $f\in C_b(\erre^n)$, $f\geq 0$, and by density to
  all bounded Borel $f$ with $f \geq 0$ $\mu$-a.e..
\end{proof}

The operator $L_0$ is $\omega$-dissipative in $L^2(\erre^n,\mu)$. More
precisely, we have
$$
\int_{\erre^n} (L_0f)f\,d\mu \leq 
-\frac12\int_{\erre^n}|(\sigma\sigma^*)^{1/2}Df|^2\,d\mu
+ \omega \int_{\erre^n} f^2\,d\mu
\qquad
\forall f \in C^2_b(\erre^n),
$$
as follows by (\ref{eq:ex}) and $L_0(f^2) = 2(L_0f)f +
|(\sigma\sigma^*)^{1/2}Df|^2$.
\par\noindent
Note also that for each $h\in C^2_b(\erre^n)$ the function
\begin{equation}
\label{eq:res}
\varphi(x) = \E_{0,x}\int_0^\infty e^{-\lambda t}h(X(t))\,dt  
\end{equation}
is in $C^2_b(\erre^n)$ and satisfies the equation
$$
\lambda \varphi - L_0\varphi = h
$$
in $\erre^n$. Hence $R(\lambda I - L_0)$ is dense in
$L^2(\erre^n,\mu)$ and since $L_0$ is closable, its closure
$L:=\overline{L_0}$ is $\omega$-$m$-dissipative, i.e. $-\omega I+L$ is
$m$-dissipative. Since, by (\ref{eq:res}), $(\lambda I - L)^{-1}$ is
the resolvent of the infinitesimal generator of $P_t$, we also infer
that $L$ is just the infinitesimal generator of $P_t$. We have
thus proved the following result.
\begin{lemma}
The infinitesimal generator of $P_t$ in $L^2(\erre^n,\mu)$ is
$L$. Moreover one has
$$
\int_{\erre^n}(Lf)f\,d\mu \leq 
- \frac12 \int_{\erre^n}|\sigma^*Df|^2\,d\mu
+ \omega \int_{\erre^n} f^2\,d\mu
$$
for all $f \in L^2(\erre^n,\mu)$.
\end{lemma}
Taking into account that $L$ is the closure (i.e.
Friedrichs' extension) of $L_0$ in $L^2(\erre,\mu)$, it follows that
for each $f \in D(L)$ we have
$$
L f = \frac12 \tr[\sigma\sigma^*D^2f]
+ \ip{rx}{Df}_{\erre^n}
$$
in $\mathcal{D}'(\erre^n)$, where $Df$, $D^2f$ are taken in the sense
of distributions. In particular, it follows by the previous lemma that 
$$
(\sigma\sigma^*)^{1/2}f \in W^{1,2}(\erre^n,\mu),
\qquad
f \in W^{2,2}_{loc}(\Xi)
$$
for each $f \in D(\overline{L})$, where
$\Xi=\{x\in\erre^n:\;\tr[\sigma\sigma^*](x)>0\}$.

\medskip

We are now going to apply Theorem \ref{thm:abs} to the operator
$N=-L+rI$ on the set
$$
\K_g = \Big\{ \varphi\in L^2(\erre^n,\mu): \; \varphi(x) \geq g(x)
\;\; \mu\mathrm{-a.e.}\Big\}.
$$
The function $g:\erre^n\to\erre$ is assumed to satisfy the following
conditions:
\begin{equation}
\label{eq:cond-g}
Dg \in L^\infty(\erre^n), \qquad
\tr[\sigma\sigma^*D^2g] \in \mathcal{M}(\erre^n), \qquad
\tr[\sigma\sigma^*D^2g] \geq 0 \;\;\mathrm{in\ } \mathcal{M}(\erre^n),
\end{equation}
where $\mathcal{M}(\erre^n)$ is the space of bounded Radon measures on
$\erre^n$.

Payoff functions that can be covered in this setting include so-called
Margrabe options (with payoff $g(x)=(x_i-\lambda x_j)^+$, for given
$\lambda>0$ and $i\neq j \leq n$) and basket put options. We shall
focus, as an illustration of the theory, on the latter case, for which
$$
g(x_1,\ldots,x_n) = \Big(k-\sum_{j=1}^n \lambda_jx_j\Big)^+, \qquad
\sum_{j=1}^n \lambda_j=1.
$$
In this case the first two conditions in (\ref{eq:cond-g}) are
obviously satisfied and
$$
\tr[\sigma\sigma^*D^2g] =
\Big(\sum_{i,j=1}^n a_{ij}\lambda_i\lambda_j\Big) \delta \geq 0
$$
in $\mathcal{M}(\erre^n)$, where $a=\sigma\sigma^*$ and $\delta$ is
the Dirac measure. Moreover,
$$
D_ig(x) = -\lambda_i H\Big(r - \sum_{j=1}^n \lambda_jx_j\Big),
\qquad i=1,\ldots,n,
$$
where $H$ is the Heaviside function, i.e. $H(r)=1$ for $r\geq 0$ and
$H(r)=0$ otherwise.
The operator $N$ has a natural extension to functions $g$ satisfying
the first two conditions in (\ref{eq:cond-g}) through the formula
$$
(Ng)\varphi = \int_{\erre^n} gN^*\varphi\,d\mu
\qquad \forall \varphi \in D(N^*).
$$
In our case one has
$$
Ng = -\frac12 \Big(\sum_{i,j=1}^n a_{ij}\lambda_i\lambda_j\Big) \delta
- r\sum_{j=1}^n \lambda_jx_j \, H\Big(k - \sum_{j=1}^n \lambda_jx_j\Big)
+ rg.
$$

\begin{lemma}
  Assume that (\ref{eq:hyp-sigma}) and (\ref{eq:cond-g}) are verified.
  Then the operator $N+\mathcal{N}_g$ is $\omega$-$m$-accretive in
  $L^2(\erre^n,\mu)$.
\end{lemma}
\begin{proof}
We only have to prove that $|(N_\lambda g)^+|_{L^2(\erre^n,\mu)}$ is
bounded for all $\lambda \in (0,\omega^{-1})$, as required by
Theorem \ref{thm:abs}. Set $g_\lambda=(I+\lambda N)^{-1}g$, i.e.
\begin{equation}
\label{eq:Ng}
g_\lambda + \lambda N g_\lambda = g, \qquad N_\lambda g = N g_\lambda.  
\end{equation}
Then we have
$$
(1+\lambda r)g_\lambda(x)
- \frac{\lambda}{2}\sum_{i,j=1}^n a_{ij} D^2_{ij}g_\lambda(x)
- r\sum_{i=1}^n x_iD_ig_\lambda(x) = g(x)
$$
in $\mathcal{D}'(\erre^n)$.
As seen earlier, $Ng=-\frac12 \tr[\sigma\sigma^*D^2g]-\ip{rx}{Dg}+rg$ in
$\mathcal{D}'(\erre^n)$ and by assumption (\ref{eq:cond-g}) we have
that $(Ng)^+ = (-r\ip{x}{Dg}+rg)^+$ (where $\nu^+$ denotes the positive
part of the measure $\nu$). Since $Dg\in L^\infty(\erre^n,dx)$ we
conclude that $(Ng)^+ \in L^2(\erre^n,\mu)$.
\par\noindent
Approximating $g$ by a sequence $g_\varepsilon\in D(N)$ we may assume
that $g\in D(N)$ and also $Ng_\lambda \in D(N)$. We set $\psi_\lambda
= Ng_\lambda$ and so (\ref{eq:Ng}) yields
$$
(1+\lambda r)\psi_\lambda(x) 
- \frac{\lambda}{2}\sum_{i,j=1}^n a_{ij}(x) D^2_{ij} \psi_\lambda(x)
-r \sum_{i=1}^n x_i D_i\psi_\lambda(x) = Ng(x).
$$
Let us set $\psi_\lambda = \psi^1_\lambda + \psi^2_\lambda$, with
\begin{eqnarray*}
(1+\lambda r)\psi^1_\lambda(x) 
- \frac{\lambda}{2} \sum_{i,j=1}^n a_{ij} D^2_{ij} \psi^1_\lambda
-r \sum_{i=1}^n x_i D_i\psi^1_\lambda &=& (Ng)^+ \\
(1+\lambda r)\psi^2_\lambda(x) 
- \frac{\lambda}{2} \sum_{i,j=1}^n a_{ij} D^2_{ij} \psi^2_\lambda
-r \sum_{i=1}^n x_i D_i\psi^2_\lambda &=& (Ng)^-,
\end{eqnarray*}
where the first equation is taken in $L^2(\erre^n,\mu)$ and the second
in $\mathcal{D}'(\erre^n)$.  By the maximum principle for elliptic
equations we infer that $\psi_\lambda^1\geq 0$, $\psi_\lambda^2\geq
0$, hence $\psi_\lambda^1 = \psi_\lambda^+$ and
$\psi_\lambda^2=\psi_\lambda^-$.
This implies that $\psi_\lambda^+=(N_\lambda g)^+$ is the solution
$\psi_\lambda^1$ of
$$
\psi_\lambda^1 + \lambda N\psi_\lambda^1 = (Ng)^+.
$$
But the solution of this equation satisfies
$$
|\psi_\lambda^1|^2_{L^2(\erre^n,\mu)} \leq
\frac{|(Ng)^+|_{L^2(\erre^n,\mu)}}{1-\lambda\omega},
$$
hence $\{(N_\lambda g)^+\}_\lambda$ is bounded as claimed.
\end{proof}

Applying Corollary \ref{cor:abs} we obtain the following existence
result for the value function of the optimal stopping problem, i.e.
for the price of the American option.
\begin{coroll}\label{cor:craxo}
  Assume that conditions (\ref{eq:hyp-sigma}), (\ref{eq:cond-g}) hold
  and that $g\in \overline{D(N)}=L^2(\erre^n,\mu)$. Then the backward
  variational inequality associated to the optimal stopping problem
  (\ref{eq:aop}), i.e.
\begin{equation}
\label{eq:caldo}
\left\{
\begin{array}{l}
\ds \frac{\partial u}{\partial t} - Nu - \mathcal{N}_g(u) \ni 0,
\quad \mathrm{a.e.\ } t\in (0,T),\\[8pt]
u(T)=g,
\end{array}\right.
\end{equation}
admits a unique generalized (mild) solution $u$ in
$C([0,T];L^2(\erre^n,\mu))$.  Moreover, if $g\in D(N)$, then $u\in
W^{1,\infty}([0,T];L^2(\erre^n,\mu))$ is the unique strong solution of
(\ref{eq:caldo}). Furthermore, if the law of the solution of
(\ref{eq:sde}) is absolutely continuous with respect to $\mu$, then
the value function $v$ coincides with $u$ for all $s\in[t,T]$ and
$\mu$-a.e. $x\in\erre^n$.
\end{coroll}
Let us remark that the last assertion of the corollary is included for
completeness only, as we do not know of any option whose payoff $g$ is
smooth enough so that $g\in D(N)$.  On the other hand, the general
case $g\in \overline{D(N)}$ covered in the corollary happens for
virtually all payoff functions $g$. Then the solution is just the
limit of the following backward finite difference scheme:
$$
v_i=\theta_{M-i}, \quad
\theta_{i+1}+hN\theta_{i+1} + \mathcal{N}_g(\theta_{i+1}) \ni \theta_i, \quad
\theta_0=v_0, \quad
h={T/M}.
$$
This discretized elliptic variational inequality can be solved via the
penalization scheme
$$
\theta^\varepsilon_{i+1} + h N\theta^\varepsilon_{i+1}
-\frac1\varepsilon (\theta_{i+1}^\varepsilon-g)^- = \theta^\varepsilon_i,
\qquad i=0,1,\ldots,M-1,
$$
or via the bounded penalization scheme (see e.g. \cite{BrNi})
$$
\theta^\varepsilon_{i+1} + h N\theta^\varepsilon_{i+1}
+g_1 \frac{\theta^\varepsilon_{i+1}-g}
          {\varepsilon+|\theta^\varepsilon_{i+1}-g|} = 
\theta^\varepsilon_i + g_1,
\qquad i=0,1,\ldots,M-1,
$$
where $g_1$ is an arbitrary parameter function. Therefore the
characterization of the option price given by Corollary
\ref{cor:craxo} is also constructive, that is, it is
guaranteed to be the unique limit of very natural finite difference
approximation schemes, that can be implemented numerically. A
completely analogous remark applies also to the cases treated in the
next subsections.

\subsection{American options on assets with stochastic volatility}
Consider the following model of asset price dynamics with stochastic
volatility under a risk neutral measure $\mathbb{Q}$:
\begin{eqnarray*}
dX(t) &=& \sqrt{V(t)} X(t)\,dW_1(t) \\
dV(t) &=& \kappa(\theta-V(t))\,dt + \eta\sqrt{V(t)}\,dW_2(t),
\end{eqnarray*}
where $W(t)=(W_1(t),W_2(t))$ is a 2-dimensional Wiener process with
identity covariance matrix (the more general case of correlated Wiener
processes is completely analogous), $\kappa$, $\theta$, $\eta$ are
positive constants, and the risk-free interest rate is assumed to be
zero. Moreover, in order to ensure that $V(t) \geq 0$
$\mathbb{Q}$-a.s. for all $t\in[0,T]$, we assume that
$2\kappa\theta>\eta^2$ (see e.g.  \cite{Hes93}).

It is convenient to use the transformation $x(t)=\log X(t)$, after
which we can write (by a simple application of It\^o's lemma)
$$
dx(t) = -V(t)/2\,dt + \sqrt{V(t)}\,dW_1(t).
$$
Define $Y(t)=(x(t),V(t))$. Then we have
\begin{equation}
\label{eq:Yp}
dY(t) = A(Y(t)) + G(Y(t))\,dW(t),  
\end{equation}
where $A: \erre^2 \ni (x,v) \mapsto
(-v/2,\kappa(\theta-v))\in\erre^2$ and $G:\erre^2 \ni (x,v) \mapsto
\mathrm{diag}(\sqrt{v},\eta\sqrt{v})\in L(\erre^2,\erre^2)$. The price
of an American contingent claim on $X$ with payoff function
$g:\erre\to\erre$ is the value function $v$ of an optimal stopping
problem, namely
\begin{equation}
  \label{eq:prezzo}
v(t,x,v) = \sup_{\tau\in\mathfrak{M}} \E_{t,(x,v)}[\tilde{g}(Y(\tau))],  
\end{equation}
where $\tilde{g}(x,v)\equiv g(e^x)$ and $\mathfrak{M}$ is the set of
all stopping times $\tau$ such that $\tau\in [s,T]$ $\mathbb{Q}$-a.s..

The Kolmogorov operator $L_0$ associated to (\ref{eq:Yp}) is given by
$$
L_0 f = \frac12 v f_{xx} + \frac12 \eta^2 v f_{vv} - \frac12v f_x +
\kappa(\theta-v)f_v, \qquad f \in C^2_b(\erre^2),
$$
and its adjoint $L_0^*$ takes the form
\begin{equation}
\label{eq:pippo}
L_0^*\rho = \frac12 v \rho_{xx} + \frac12 \eta^2(v\rho)_{vv} + \frac12v\rho_x
- \kappa((\theta-v)\rho)_v, \qquad \rho \in C^2_b(\erre^2),
\end{equation}
Following the same strategy as above, we look for an excessive measure
of the form
$$
\mu(dx,dv) = a\rho(x,v)\,dx\,dv,
\qquad
\rho(x,v) = \frac1{1+x^2+v^2},
$$
where $a^{-1}=\int_{\erre^2} \rho(x,v)\,dx\,dv$.

Some calculations involving (\ref{eq:pippo}) reveal that
$$
\sup_{(x,v)\in\erre\times\erre_+} \frac{L_0^*\rho(x,v)}{\rho(x,v)} =
\omega < \infty,
$$
i.e. $\mu$ is an infinitesimally excessive measure for $L_0$ on
$\Xi=\erre\times\erre_+$.  Then the transition semigroup
$$
P_tf(x,v) = \E_{0,(x,v)}f(x(t),V(t)), \qquad f\in C^2_b(\Xi),
$$
extends by continuity to $L^2(\Xi,\mu)$, and the operator $L_0$ with
domain $C^2_b(\Xi))$ is $\omega$-dissipative in $L^2(\Xi,\mu)$.
Arguing as above (see Lemma \ref{lem:clos}), the closure
$L$ of $L_0$ is $\omega$-$m$-dissipative in $L^2(\Xi,\mu)$
and
$$
\int_\Xi (Lf) f\,d\mu \leq 
- \frac{\eta^2+1}{2} \int_\Xi v(f_x^2+f_v^2)\,d\mu
+ \omega \int_\Xi f^2\,d\mu.
$$
The operator $N=-L$ is therefore $\omega$-$m$-accretive and
formally one has
\begin{eqnarray}
N\tilde{g} &=& -\frac12 v(e^{2x}g''(e^x)+e^xg'(e^x))
               +\frac12 v e^xg'(e^x) \nonumber\\
&=& -\frac12 v e^{2x} g''(e^x).
\label{eq:pluto}
\end{eqnarray}
The previous expression is of course rigorous if $g$ is smooth and
$N\tilde{g}\in L^2(\Xi,\mu)$, but in general (i.e. for $\tilde{g} \in
L^2(\Xi,\mu)$) is has to be interpreted in the sense of distributions
on $\Xi$ in order to be meaningful.
\par\noindent
We shall assume that the payoff function $g$ is convex on $\erre$,
more precisely,
\begin{equation}
  \label{eq:pape}
g''\in\mathcal{M}(\erre), \qquad g''\geq 0,
\end{equation}
where $\mathcal{M}(\erre)$ is the space of finite measures on $\erre$.
Note that the typical payoff of a put or call option is covered by
these assumptions.  Equation (\ref{eq:pluto}) implies that
$N\tilde{g}\in\mathcal{D}'(\erre)$ and $N\tilde{g} \leq 0$ in
$\mathcal{D}'(\erre)$, hence $N\tilde{g}$ is a negative measure and so
the hypotheses of Theorem \ref{thm:abs} are met.
Thus, defining
$\K_g=\{\varphi\in L^2(\Xi,\mu):\;\; \varphi\geq
\tilde{g}\;\;\mu\textrm{-a.e.}\}$, it follow that the operator
$N+\mathcal{N}_g$ is $\omega$-$m$-accretive on $H=L^2(\Xi,\mu)$. This
yields
\begin{coroll}
  Assume that (\ref{eq:pape}) holds. Then the backward variational
  inequality
  \begin{equation}
  \label{eq:rino}
  \frac{\partial u}{\partial t} - Nu - \mathcal{N}_g(u) \ni 0
  \end{equation}
  on $H_T=[0,T]\times L^2(\Xi,\mu)$ with terminal condition
  $u(T)=\tilde{g}$ has a unique generalized (mild) solution $u\in
  C([0,T],L^2(\Xi,\mu))$.  Moreover, if $g\in D(N)$, then
  (\ref{eq:rino}) has a unique strong solution $u\in
  W^{1,\infty}([0,T],L^2(\Xi,\mu))$.  Furthermore, if the law of the
  solution of (\ref{eq:Yp}) is absolutely continuous with respect to
  $\mu$, then the value function $v$ defined in (\ref{eq:prezzo})
  coincides with $u$ for all $s\in[t,T]$ and $\mu$-a.e. $(x,v)\in\Xi$.
\end{coroll}

\subsection{Asian options with American feature}
Let the price process $X$ of a given asset satisfy the following
stochastic differential equation, under an equivalent martingale
measure $\mathbb{Q}$:
$$
dX = rX\,dt + \sigma(X)\,dW(t), \qquad X(0)=x.
$$
Here we consider the problem of pricing a ``regularized'' Asian
options with American feature, that is we look for the value function
$v$ of the optimal stopping problem
\begin{equation}
\label{eq:ringhio}
v(x) = \sup_{\tau\in\mathfrak{M}} \E_x
\left(k - \frac{1}{\tau+\delta} \int_0^\tau X_s\,ds\right)^+,
\end{equation}
where $k\geq 0$ is the strike price, $\delta>0$ is a ``small''
regularizing term, $\mathfrak{M}$ is the set of stopping times between
$0$ and $T$, and $\E_x$ stands for expectation w.r.t. $\mathbb{Q}$,
conditional on $X(0)=x$. The standard Asian payoff corresponds to
$\delta=0$.  Unfortunately we are not able to treat with our methods
this limiting situation, as it gives rise to a singularity in the
obstacle function of the associated variational inequality, or, in the
approach we shall follow here, in the Kolmogorov operator of an
associated stochastic system. However, it is clear that for small
values of $\delta$ the value function $v$ in (\ref{eq:ringhio}) is a
good approximation of the option price, at least for optimal exercise
times that are not of the same order of magnitude of $\delta$.

Let us define the auxiliary processes
$$
Y(t) = \frac{1}{t+\delta}\int_0^t X(s)\,ds
$$
and $S(t)=t$. Then we have
$$
\left\{
\begin{array}{l}
dX(t) = rX(t)\,dt + \sigma(X(t))\,dW(t)\\[6pt]
\ds dY(t) = \frac{X(t)-Y(t)}{S(t)+\delta}\,dt\\[6pt]
dS(t)=dt
\end{array}
\right.
$$
with initial conditions $X(0)=x$, $Y(0)=0$, $s(0)=0$.  This system can
be equivalently written in terms of the vector $Z=(X,Y,S)$ as
\begin{equation}
\label{eq:zizou}
dZ(t) = A(Z(t))\,dt + G(Z(t))\,dW(t),
\qquad Z(0)=(x,0,0),
\end{equation}
where $A:\erre^3\to\erre^3$, $A:(x,y,s)\mapsto
(rx,(s+\delta)^{-1}(x-y),1)$ and $G:\erre^3\to L(\erre,\erre^3) \simeq
\erre^3$, $G(x,y,s)=(\sigma(x),0,0)$. Therefore (\ref{eq:ringhio}) is
equivalent to
$$
v(x) = \sup_{\tau\in\mathfrak{M}} \E_x g(Z(\tau)),
$$
where $g:(x,y,s)\to (k-y)^+$ and $\E_x$ stands for $\E_{(x,0,0)}$.

As in the previous cases, we shall look for an excessive measure of
$L_0$, the Kolmogorov operator associated to (\ref{eq:zizou}), which
is given by
$$
L_0f = \frac12 \sigma^2(x)D^2_{xx}f + rxD_xf
+ \frac{x-y}{s+\delta}D_yf + D_sf,
\qquad f\in C^2_b(\erre^3).
$$
Then the adjoint of $L_0$ can be formally written as
$$
L_0^*\rho = \frac12 D^2_{xx}(\sigma^2(x)\rho)
- r D_x(x\rho) - D_y\Big( \frac{x-y}{s+\delta}\rho \Big) - D_s\rho.
$$
In analogy to previous cases, some calculations reveal that, under the
assumptions (\ref{eq:hyp-sigma}) on $\sigma$, there exists a measure
$\mu$ of the type $\mu(dx,dy,ds)=\rho(x,y,s)\,dx\,dy\,ds$,
$$
\rho(x,y,s) = \frac{a}{(1+|x|)^{2n+1}(1+|x-y|)^{2n+1}(1+s)^2},
\qquad a^{-1} = \int_{\erre^3} \rho(z)\,dz,
$$
such that $L_0^*\rho \leq \omega\rho$ for some $\omega\in\erre$.
Arguing as before, we conclude that $\mu$ is an excessive measure for
the semigroup $P_t$ generated by the stochastic equation
(\ref{eq:zizou}), and that $L$, the closure of $L_0$ in $L^2(\erre^3,\mu)$,
is the infinitesimal generator of $P_t$.

We are now in the setting of section \ref{sec:viosp}, i.e. we can
characterize the option price as the (generalized) solution of a
suitable variational inequality. Details are left to the reader.

\subsection{Path-dependent American options}
We shall consider a situation where the price dynamics is
non-Markovian as it may depend on its history, and the payoff function
itself is allowed to depend on past prices. We should remark, however,
that in the present setup we still cannot cover Asian options of the
type discussed in the previous subsection, with $\delta=0$.

Consider the following price evolution of $n$ assets under a
risk-neutral measure $\mathbb{Q}$:
\begin{equation}
\label{eq:delay}
\left\{
\begin{array}{l}
dX(t)=rX(t) + \sigma(X(t),X_s(t))\,dW(t),\quad 0 \leq t \leq T \\
X(0)=x_0, \qquad X_s(0)=x_1(s), \quad -T\leq s \leq 0,
\end{array}
\right.
\end{equation}
where $X_s(t)=X(t+s)$, $s\in(-T,0)$, $W$ is a standard Wiener process
on $\erre^n$ and $\sigma:\erre^n\times L^2([-T,0],\erre^n) \to
L(\erre^n,\erre^n)$.  Let us consider an American contingent claim
with payoff of the type $g:\erre^n\times L^2([-T,0]\to\erre$, whose
price is equal to the value function of the optimal stopping problem
$$
v(s,x_0,x_1) = \sup_{\tau\in\mathfrak{M}}
\E_{s,(x_0,x_1)} [e^{-r\tau}g(X(\tau),X_s(\tau))],
$$
where the notation is completely analogous to the previous subsection.
An example that can be covered by this functional setting is
$g(x_0,x_1)=\alpha_0g_0(x_0)+\alpha_2g_1(x_1)$, with $\alpha_1$,
$\alpha_2\geq 0$ and $g_0(x_0)=(k_0-x_0)^+$ and
$g_1(x_1)=(k_1-\int_{-T}^0 x_1(s)\,ds)^+$.

Let us now rewrite (\ref{eq:delay}) as an infinite dimensional
stochastic differential equation on the space $H=\erre^n \times
L^2([-T,0],\erre^n)$.
Define the operator $A:D(A)\subset H \to H$ as follows:
\begin{eqnarray*}
A: (x_0,x_1) &\mapsto& (rx_0,x_1') \\
D(A) &=& \{ (x_0,x_1)\in H;\; x_1 \in W^{1,2}((-T,0),\erre^n),\; x_1(0)=x_0 \}.
\end{eqnarray*}
Setting $G(x_0,x_1)=(\sigma(x_0,x_1),0)$, let us consider the stochastic
differential equation on $H$
\begin{equation}
\label{eq:delay2}
dY(t) = AY(t)\,dt + G(Y(t))\,dW(t)
\end{equation}
with initial condition $Y(0)=(x_0,x_1)$.  The evolution equation
(\ref{eq:delay2}) is equivalent to (\ref{eq:delay}) in the following
sense (see \cite{choj78}): if $X$ is the unique solution of
(\ref{eq:delay}), then $Y(t)=(X(t),X_s(t)$ is the solution of
(\ref{eq:delay2}). Note that (\ref{eq:delay2}) has a unique solution
if $G$ is Lipschitz on $H$.
Finally, regarding $g$ as a real-valued function defined on $H$, we
are led to study the optimal stopping problem in the Hilbert space $H$
$$
v(s,x) = \sup_{\tau\in\mathfrak{M}} \E_{s,x}[e^{-r\tau}g(Y(\tau))].
$$
The Kolmogorov operator $L_0$ associated to (\ref{eq:delay2}) has the form,
on $C^2_b(H)$,
\begin{eqnarray*}
L_0\varphi(x_0,x_1) &=&
\frac12 \tr[\sigma\sigma^*(x_0,x_1)D^2_{x_0}\varphi(x_0,x_1)]
+ \ip{rx_0}{D_{x_0}\varphi(x_0,x_1)}_{\erre^n} \\
&& + \int_{-T}^0 \ip{x_1(s)}{D_{x_1}\varphi(x_0,x_1(s))}_{\erre^n}\,ds.
\end{eqnarray*}
We look for an excessive measure $\mu$ for $L_0$ of the form
$\mu=\nu_1\otimes\nu_2$, where $\nu_1$, $\nu_2$ are probability
measures on $\erre^n$ and $L^2([-T,0],\erre^n)$, respectively. In
particular, we choose
$$
\nu_1(dx_0) = \rho(x_0)\,dx_0,
\qquad
\rho(x_0) = \frac{a}{1+|x_0|^{2n}},
\qquad
a = \left( \int_{\erre^n} \frac{1}{1+|x_0|^{2n}}\,dx_0 \right)^{-1},
$$
and $\nu_2$ a Gaussian measure on $L^2([-T,0],\erre^n)$.
Setting $H_0=L^2([-T,0],\erre^n)$, a simple calculation reveals that
\begin{eqnarray}
\int_H L_0\varphi\,d\mu &=&
\frac12 \int_{H_0}d\nu_2\int_{\erre^n}
            \tr[\sigma\sigma^*D^2_{x_0}\varphi]\,d\nu_1
+ r \int_{H_0}d\nu_2\int_{\erre^n}
          \ip{x_0}{D_{x_0}\varphi}_{\erre^n}\,d\nu_1\nonumber\\
&& + \int_{\erre^n}d\nu_1 \int_{H_0}
          \ip{x_1}{D_{x_1}\varphi}_{H_0}\,d\nu_2\nonumber\\
&=& \frac12 \int_{H_0}d\nu_2 \int_{\erre^n}
                  \varphi D^2_{x_0}(\sigma\sigma^*\rho)\,dx_0
- r \int_{H_0}d\nu_2 \int_{\erre^n} \varphi D_{x_0}(x_0\rho)\,dx_0\nonumber\\
&& + \int_{\erre^n} d\nu_1 \int_{H_0} \ip{x_1}{D_{x_1}\varphi}_{H_0}\,d\nu_2.
\label{eq:fidanken}
\end{eqnarray}
We shall assume that
\begin{equation}
\label{eq:hass}
\begin{array}{l}
\;\;\sigma \in C^2(\erre^n\times L^2([-T,0],\erre^n)) \cap
\mathrm{Lip}(\erre^n\times L^2([-T,0],\erre^n)), \\[4pt]
\begin{array}{ll}
\sigma(x_0,x_1) \leq C(|x_0|+|x_1|_{H_0}) & \forall (x_0,x_1)\in H, \\[4pt]
|\sigma_{x_i}(x_0,x_1)| + |\sigma_{x_ix_j}(x_0,x_1)| \leq C &
\forall (x_0,x_1)\in H, \;\; i,j=1,2.
\end{array}
\end{array}
\end{equation}
Note that these conditions also imply existence and uniqueness of a
solution for (\ref{eq:delay}).
\par\noindent
Taking into account that $\int_{H_0} |x_1|^{2m}\,d\nu_2<\infty$ and that
$$
\int_{H_0} \ip{x_1}{D_{x_1}\varphi}_{H_0}\,d\nu_2 = 
- \int_{H_0} \varphi(n-\ip{Q^{-1}x_1}{x_1}_{H_0})\,d\nu_2
$$
(where $Q$ is the covariance operator of $\nu_2$), we see by
(\ref{eq:fidanken}) and (\ref{eq:hass}) that there exists $\omega\geq
0$ such that
\begin{equation}
  \label{eq:boss}
  \int_H L_0\varphi\,d\mu \leq \omega \int_H \varphi\,d\mu
\end{equation}
for all $\varphi\in C^2_b(H)$, $\varphi\geq 0$.

The operator $L_0$ is thus closable and $\omega$-dissipative in $L^2(H,\mu)$.
Moreover, (\ref{eq:boss}) implies that
\begin{equation}
  \label{eq:hogg}
\int_H (L_0\varphi)\varphi\,d\mu \leq
- \frac12 \int_H |(\sigma\sigma^*)^{1/2}D_{x_0}\varphi|^2\,d\mu
+ \omega \int_H \varphi^2\,d\mu
\quad \forall \varphi \in C^2_b(H).
\end{equation}
Since one has, for $\lambda>\omega$,
$$
(\lambda I - L_0)^{-1}\varphi = 
\E\int_0^\infty e^{-\lambda t}\varphi(X(t),X_s(t))\,dt
\quad \forall \varphi \in C^2_b(H),
$$
we infer that $R(\lambda I - L_0)$ is dense in $L^2(H,\mu)$ and so the
closure $L$ of $L_0$ is $\omega$-$m$-dissipative in
$L^2(H,\mu)$, and it is the infinitesimal generator of the transition
semigroup $P_t$ defined by (\ref{eq:delay2}). We set $N=-L+rI$.

Furthermore, let us assume that
\begin{equation}
  \label{eq:james}
  g(\cdot,x_1) \in \mathrm{Lip}(\erre^n),
  \qquad
  D^2_{x_0}g(\cdot,x_1) \in \mathcal{M}(\erre^n) \quad \forall x_1 \in H_0,
\end{equation}
\begin{equation}
  \label{eq:bond}
  \tr[\sigma\sigma^* D^2_{x_0}g](\cdot,x_1) \geq 0 \quad \forall x_1 \in H_0,
\end{equation}
where (\ref{eq:bond}) is taken in the sense of distributions (or
equivalently in the sense of $\mathcal{M}(\erre^n)$).
This implies, as in previous cases, that condition (\ref{eq:gt2}) is
satisfied.
\par\noindent
In particular, note that (\ref{eq:james}) and (\ref{eq:bond}) hold if
$g=\alpha_0g_0+\alpha_1g_1$, as in the example mentioned above.
Assumptions (\ref{eq:james}) and (\ref{eq:bond}) imply that
\begin{eqnarray*}
Ng &=& -\frac12 \tr[\sigma\sigma^*D^2g] - r\ip{x_0}{D_{x_0}g}_{\erre^n}
- \ip{x_1}{D_{x_1}g}_{H_0} + rg\\
&\leq& - r\ip{x_0}{D_{x_0}g}_{\erre^n}
- \ip{x_1}{D_{x_1}g}_{H_0} + rg,
\end{eqnarray*}
hence $(Ng)^+\in L^2(H,\mu)$, because $\ip{x_0}{D_{x_0}g}_{\erre^n}$,
$\ip{x_1}{D_{x_1}g}_{H_0}$, $g \in L^2(H,\mu)$.

Once again the results established in sections \ref{sec:viosp} allow us to
characterize the price of the American option as solution (mild, in general,
as the typical payoff function $g$ is not smooth) of the backward variational
inequality on $[0,T]$
$$
{d\varphi \over dt} - N\varphi - \mathcal{N}_g(\varphi) \ni 0
$$
with terminal condition $\varphi(T)=g$, where
$\mathcal{N}_g$ is the normal cone to
$$
\mathcal{K}_g = \{x \in H:\; \varphi(x) \geq g(x)\;\mu\textrm{-a.e.}\}.
$$

\subsection*{Acknowledgements}
This work done during the visit of the first author as Mercator
Gastprofessor at the Institute of Applied Mathematics, University of
Bonn. The second author gratefully acknowledges the financial support
of the SFB 611, Bonn, of IMPAN, Warsaw and IH\'ES, Bures-sur-Yvette
through an IPDE fellowship, and of the ESF through grant AMaMeF 969.


\def\polhk#1{\setbox0=\hbox{#1}{\ooalign{\hidewidth
  \lower1.5ex\hbox{`}\hidewidth\crcr\unhbox0}}}

\end{document}